\documentclass[12pt]{amsart}
\usepackage{latexsym}
\usepackage{amssymb,times,fullpage}
\newtheorem{thm}{Theorem}
\newtheorem{prop}[thm]{Proposition}

\newtheorem{cor}[thm]{Corollary}

\theoremstyle{remark}
\newtheorem{rem}[thm]{Remark}
\newtheorem{ex}[thm]{Example}

\theoremstyle{definition}
\newtheorem{defn}[thm]{Definition}

\newcommand{\R}{\mathbb{ R}}

\newcommand{\Z}{\mathbb{ Z}}

\title[Quasi-homomorphisms and stable lengths in mapping class groups]{Quasi-homomorphisms and stable lengths in\\ mapping class groups}
\author{D.~Kotschick}
\address{Mathematisches Institut, Ludwig-Maximilians-Universit\"at M\"unchen,
Theresienstr.~39, 80333 M\"unchen, Germany}
\email{dieter@member.ams.org}
\thanks{The author is a member of the 
{\sl European Differential Geometry Endeavour} (EDGE), Research 
Training Network HPRN-CT-2000-00101, supported by The European Human 
Potential Programme}
\date{August 26, 2003; MSC 2000: primary 20F69, secondary 20F12, 57M07}

\begin{document}

\begin{abstract}
    We give elementary applications of quasi-homomorphisms to 
    growth problems in groups. A particular case concerns the 
    number of torsion elements required to factorise a given element 
    in the mapping class group of a surface.
\end{abstract}

\maketitle

\section{Introduction}

Let $G$ be a group, and $S\subset G$ an arbitrary subset. The 
$S$-length $l_{S}(g)$ of $g\in G$ is the minimal number of factors in 
a factorisation\footnote{The length is infinite if no such 
factorisation exists.} of $g$ into elements of $S$. For a 
homomorphism $\varphi\colon G\rightarrow\R$ set 
$C(\varphi,S)=\sup_{s\in S}\vert\varphi(s)\vert$. Then 
\begin{equation}\label{e:length}
    n\cdot\vert\varphi(g)\vert = \vert\varphi(g^{n})\vert\leq 
    l_{S}(g^{n})\cdot C(\varphi,S) \ .
    \end{equation}
If $\varphi(g)\neq 0$ and $\varphi\vert_{S}$ is bounded\footnote{for 
example if $S$ is finite, or, more generally, falls into finitely 
many conjugacy classes}, then we obtain a positive lower bound for the 
$S$-length $l_{S}(g^{n})$ of $g^{n}$, which is linear in $n$. Hence we 
obtain a positive lower bound for the stable $S$-length $\vert\vert 
g\vert\vert_{S}=\lim_{n\rightarrow\infty}\frac{l_{S}(g^{n})}{n}$ of $g$:
\begin{equation}\label{e:slength}
    \vert\vert g \vert\vert_{S}
    \geq\frac{\vert\varphi(g)\vert}{C(\varphi,S)} \ .
    \end{equation}
In many important situations there are no homomorphisms from $G$ to $\R$ 
to which one can apply this argument, either because there are no 
homomorphisms at all, or they all vanish on a given $g$. However, 
quasi-homomorphisms can be used in the same way as homomorphisms to 
find lower bounds on the length. If $\varphi$ is a homogeneous quasi-homomorphism 
with defect $D(\varphi)$ (see Section~\ref{s:qmorph} below for the 
definitions), then instead of~\eqref{e:length} we have 
\begin{equation}\label{e:qlength}
    n\cdot\vert\varphi(g)\vert = \vert\varphi(g^{n})\vert\leq 
    l_{S}(g^{n})\cdot C(\varphi,S) +(l_{S}(g^{n})-1)\cdot 
    D(\varphi)\ ,
    \end{equation}
which again leads to a lower bound for $l_{S}(g^{n})$ which is positive 
and linear in $n$ as soon as $\varphi(g)\neq 0$ and $\varphi\vert_{S}$ 
is bounded. For the stable $S$-length we then have\footnote{Note that 
the stable $S$-length can be defined and is finite as soon as $g$ has 
some power that is a product of elements in $S$, even if $g$ itself cannot 
be factorised in $S$.}
\begin{equation}\label{e:qslength}
   \vert\vert g\vert\vert_{S}
    \geq\frac{\vert\varphi(g)\vert}{C(\varphi,S)+D(\varphi)} \ .
    \end{equation}
This kind of estimate has been discussed extensively in the 
literature for the case when $S$ is taken to be the set of commutators 
in $G$, see Bavard~\cite{B} and the papers cited there, 
and~\cite{BF,EK,PR} for some more recent developments.

In this paper we apply the estimates~\eqref{e:qlength} and~\eqref{e:qslength}
in situations when $S$ is not necessarily the set of commutators. However, 
as the positivity of the stable commutator length is equivalent to the 
existence of homogeneous quasi-homomorphisms by Bavard's theorem~\cite{B}, 
the stable commutator length provides a kind of universal lower bound for the 
stable $S$-length for all $S$. We will see in particular that if we take $S$ to 
be the set of elements of finite order in $G$, then the corresponding stable 
torsion length is bounded below by twice the stable commutator length.

In the final section we apply our results to mapping class groups of 
surfaces. We will also exhibit elements $g$ of infinite order in 
mapping class groups for which the commutator length and the torsion length 
of $g^{n}$ are constant, and therefore the stable commutator and 
torsion lengths vanish. For these elements the stable $S$-length is positive 
for any finite $S$ by a result of Mosher~\cite{Mosher} and 
Farb--Lubotzky--Minsky~\cite{FLM}. Thus, quasi-homomorphisms give lower bounds 
for the stable $S$-length, but they only tell the whole story when $S$ is the 
set of commutators. We shall see that quasi-homomorphisms on mapping 
class groups also lead to a unified approach to several other of their 
properties, which were originally established by other means: failure 
of (weak) bounded generation~\cite{FLM}, non-arithmeticity~\cite{I} 
and superrigidity~\cite{FM,KM}.

\section{Quasi-homomorphisms}\label{s:qmorph}

A map $f\colon G\rightarrow\R$ is called a quasi-homomorphism 
if its deviation from being a homomorphism is bounded; in other words, 
there exists a constant $D(f)$, called the defect of $f$, 
such that
$$
\vert f(xy)-f(x)-f(y)\vert\leq D(f)
$$
for all $x,y\in G$. Every quasi-homomorphism can be homogenized by 
defining
$$
\varphi (g) = \lim_{n\rightarrow\infty}\frac{f(g^{n})}{n} \ .
$$
Then $\varphi$ is again a quasi-homomorphism, is homogeneous in the 
sense that $\varphi(g^{n})=n\varphi(g)$, and is constant on conjugacy 
classes. (Compare~\cite{B}, Proposition~3.3.1.) The 
inequality~\eqref{e:qlength} is immediate.

\subsection{Commutator lengths}

Consider the situation when $S=C$ is the set of commutators. 
In this case we write $c(g)$ for the commutator length $l_{C}(g)$, 
and $\vert\vert g\vert\vert_{C}$ for the stable commutator length.

For any quasi-homomorphism $\varphi$ we have
$$
\vert\varphi([x,y])\vert\leq\vert\varphi(xyx^{-1})+\varphi(y^{-1})\vert
+D(\varphi) \ .
$$
If $\varphi$ is homogeneous, then $\varphi(xyx^{-1})=\varphi(y)$ by 
constancy on conjugacy classes, and $\varphi(y^{-1})=-\varphi(y)$ by 
homogeneity. Thus $\varphi\vert_{C}$ is bounded by 
$C(\varphi,C)=D(\varphi)$, and~\eqref{e:qlength} becomes
\begin{equation}\label{e:commlength}
    n\cdot\vert\varphi(g)\vert = \vert\varphi(g^{n})\vert\leq 
    (2c(g^{n})-1)\cdot D(\varphi) \ .
    \end{equation}
In particular,
\begin{equation}\label{e:scommlength}
   \vert\vert g\vert\vert_{C}
    \geq\frac{\vert\varphi(g)\vert}{2D(\varphi)} \ .
    \end{equation}
Bavard's main theorem in~\cite{B} states that taking the supremum over the 
right hand sides for all homogeneous quasi-homomorphisms gives an 
exact calculation of the stable commutator length. In particular, if 
$\vert\vert g\vert\vert_{C}$ is positive, then there exists a 
homogeneous quasi-homomorphism $\varphi\colon G\rightarrow\R$ with 
$\varphi (g)\neq 0$.
    
\subsection{Torsion lengths}\label{ss:tor}

Now we take $S=T$, the set of torsion elements of $G$. We 
denote the torsion length $l_{T}(g)$ by $t(g)$.

Clearly every homogeneous quasi-homomorphism $\varphi$ vanishes on $T$. 
Therefore we can take $C(\varphi,T)=0$, and we find 
from~\eqref{e:qlength}:
\begin{equation}\label{e:torlength}
    n\cdot\vert\varphi(g)\vert = \vert\varphi(g^{n})\vert\leq 
    (t(g^{n})-1)\cdot 
    D(\varphi)\ .
    \end{equation}
In particular,
\begin{equation}\label{e:storlength}
   \vert\vert g\vert\vert_{T}
    \geq\frac{\vert\varphi(g)\vert}{D(\varphi)} \ .
    \end{equation}
Comparing~\eqref{e:storlength} and~\eqref{e:scommlength}, and using the 
fact that~\eqref{e:scommlength} is sharp~\cite{B}, leads to the 
following:
\begin{prop}\label{p:tor}
    Let $G$ be a group with finitely generated Abelianization. 
    For all $g\in G$ we have 
    $2\vert\vert g\vert\vert_{C}\leq\vert\vert g\vert\vert_{T}$.
    \end{prop}
\begin{proof}
    Both sides of the inequality are homogeneous, so we can replace 
    $g$ by a power if necessary.
    
    If the stable torsion length of $g$ is defined and finite, then 
    $g$ has a power that is a product of elements of finite order. 
    Therefore every homomorphism $G\rightarrow\Z$ vanishes on $g$.
    As $G$ is assumed to have finitely generated Abelianization,  
    $g$ must have a power that is in the commutator subgroup. Thus, 
    after replacing $g$ by a suitable power, we may assume that it is 
    a product of commutators and is also a product of torsion elements 
    $t_{1},\ldots,t_{k}$, with $k=t(g)$ the torsion length.
    
    For arbitrary elements $x,y$ in the commutator subgroup of $G$ one 
    has 
    $$
    \vert\vert xy\vert\vert_{C}\leq\vert\vert x\vert\vert_{C}+
    \vert\vert y\vert\vert_{C}+\frac{1}{2}
    $$
    because $(xy)^{n}$ differs from $x^{n}y^{n}$ by at most 
    $\frac{n}{2}+c$ commutators, see~\cite{B}, Proposition 3.7.1. 
    Applying this repeatedly, we obtain
    $$
    \vert\vert g\vert\vert_{C}\leq\sum_{i=1}^{k}\vert\vert 
    t_{i}\vert\vert_{C}+\frac{1}{2}(k-1) \ .
    $$
    As the stable commutator length of an element of finite order 
    vanishes, we conclude that
    $$
    t(g)=k\geq 1+2\vert\vert g\vert\vert_{C} \ .
    $$
    Applying this to $g^{n}$ instead of $g$, using the homogeneity of 
    $\vert\vert\cdot\vert\vert_{C}$, and taking the limit for 
    $n\rightarrow\infty$ we obtain the claim.
    \end{proof}

    If $g\in G$ is conjugate to some power $g^{n}$ with $n\neq 1$, then 
clearly every homogeneous quasi-homomorphism vanishes on $g$, and 
therefore its stable commutator length also vanishes. A special case 
of such an element is a product of two involutions. 
In this case we can go further and show the vanishing of the stable 
torsion length, because every power of $g$ has torsion length at 
most $2$:
\begin{ex}\label{ex:involutions}
    Suppose $g=st$, with $s$ and $t$ of order $2$. Then 
    $$
    g^{2n}=\alpha t \alpha^{-1} \cdot t
    $$
    with $\alpha = g^{n-1}s = stst\ldots sts$. As the order of an 
    element is invariant under conjugation, we 
    conclude that $t(g^{2n})\leq 2$. We also have 
    $$
    g^{2n+1}=\alpha t \alpha^{-1} \cdot tst = \alpha t \alpha^{-1} \cdot 
    tst^{-1} \ ,
    $$
    and so $t(g^{2n+1})\leq 2$ as well.
    
    We conclude that for such a product of involutions both the 
    stable torsion length and the stable commutator length vanish.
    \end{ex}
    
    \begin{ex}\label{ex:SL}
	In $SL(2,\Z)$ take 
	$$
	g=\begin{pmatrix} 2 & 1 \\ 1 & 1 \end{pmatrix} = 
	\begin{pmatrix} 0 & 1 \\ -1 & 0 \end{pmatrix} 
	    \begin{pmatrix} -1 & -1 \\ 2 & 1 \end{pmatrix} \ .
	$$ 
	Then $g$ has infinite order, but the factors on the right hand side 
	both have order $4$. Projecting to $PSL(2,\Z)$, these factors become 
	involutions, so by Example~\ref{ex:involutions} the image of any 
	power of $g$ in $PSL(2,\Z)$ has torsion length at most $2$. It 
	follows that the torsion length of $g^{n}\in SL(2,\Z)$ is at most $3$.
	\end{ex}
	In the above example $g$ is conjugate to $g^{-1}$. On the other hand,
	if an element $g\in SL(2,\Z)$ of infinite order is not conjugate to 
	its inverse, then Polterovich--Rudnick~\cite{PR} proved that there 
	is a homogeneous quasi-homomorphism which does not vanish on $g$. 
	Therefore, our discussion above shows that $g$ has positive stable 
	torsion length.
	
    In Section~\ref{s:map} we will generalize Example~\ref{ex:SL} to mapping
    class groups of surfaces of higher genus, by exhibiting elements of 
    infinite order for which the torsion lengths of all their powers are 
    equal to $2$. 
    
\subsection{Bounded generation}\label{ss:bound}

A group $G$ is said to be boundedly generated if there are finitely 
many elements $g_{1},\ldots,g_{k}$ such that any element $g\in G$ can 
be written in the form 
$
g=\prod_{i=1}^{k}g_{i}^{\alpha_{i}}
$
for some $\alpha_{i}\in\Z$. A potentially weaker property is the 
following:
\begin{defn}
    A group $G$ is weakly boundedly generated by $g_{1},\ldots,g_{k}$ 
    if for every $g\in G$ there is a number $N$ such that all powers 
    $g^{n}$ can be written in the form 
\begin{equation}\label{eq:BG}
g^{n}=\prod_{i=1}^{N}h_{i}(n)^{\alpha_{i}(n)}
\end{equation}
for some $\alpha_{i}(n)\in\Z$, and with each $h_{i}(n)$ conjugate to 
some $g_{j}\in\{g_{1},\ldots,g_{k}\}$.
\end{defn}
The point is that the $g_{j}$ be conjugated, permuted and repeated at 
will, and that $N$ may depend on $g$, though not on $n$. The 
independence of $n$ implies that the space of homogeneous 
quasi-homomorphisms on $G$ is necessarily finite-dimensional:
\begin{prop}\label{p:BG}
    If $G$ is weakly boundedly generated by $g_{1},\ldots,g_{k}$, then 
    the vector space of homogeneous quasi-homomorphisms $\varphi\colon 
    G\rightarrow\R$ is of dimension at most $k$.
    \end{prop}
\begin{proof}
    Applying $\varphi$ to~\eqref{eq:BG} and using the 
    quasi-homomorphism property $N-1$ times, we find
    $$
    \vert\varphi(g^{n})-
    \sum_{i=1}^{N}\varphi(h_{i}(n)^{\alpha_{i}(n)})\vert\leq 
    (N-1)D(\varphi) \ .
    $$
    Using the homogeneity of $\varphi$ and dividing by $n$, we obtain:
    $$
    \vert\varphi(g)-
    \sum_{i=1}^{N}\frac{\alpha_{i}(n)}{n}\varphi(h_{i}(n))\vert\leq 
    \frac{N-1}{n}D(\varphi) \ .
    $$
    Letting $n$ tend to infinity, the right hand side tends to zero 
    because $N$ is independent of $n$, and so we see that $\varphi(g)$ 
    is determined by the values of $\varphi$ on the $h_{i}(n)$, each 
    of which is, by assumption, conjugate to some $g_{j}$. As $\varphi$ is 
    constant on conjugacy classes, it is completely determined 
    by its values on $g_{1},\ldots,g_{k}$.
    \end{proof}
\begin{rem}
    The conclusion of the Proposition can be strengthened if one knows 
    that $\varphi(g_{j})$ has to vanish for all $\varphi$, for example 
    because $g_{j}$ is torsion, or is conjugate to one of its powers. 
    \end{rem}
    
\section{Mapping class groups}\label{s:map}
    
Consider the mapping class group $\Gamma_{h}$ of isotopy classes of 
orientation-preserving diffeomorphisms of a compact 
oriented surface of genus $h$. Maclachlan~\cite{Mac} proved that this 
is generated by (finitely many) torsion elements. Thus in particular its
Abelianization is finite. For $h\geq 3$ it is known that the 
Abelianization is trivial, i.~e.~the mapping class group is 
perfect~\cite{Powell}.

By using results about the Seiberg--Witten invariants of symplectic 
four-manifolds, it has been proved that the commutator length of a power of 
a Dehn twist $g\in\Gamma_{h}$ grows linearly with the exponent, see~\cite{EK} 
and the elaborations in~\cite{BK,Kork}. The argument 
admits the following slight generalization:
\begin{thm}\label{t:main}
    Let $\Gamma_{h}$ be the mapping class group of a closed oriented 
    surface $\Sigma_{h}$ of genus $h\geq 2$. If $g\in\Gamma_{h}$ is the product 
    of $k$ right-handed Dehn twists along homotopically essential 
    disjoint curves $a_{1},\ldots,a_{k}\subset\Sigma_{h}$, then 
    \begin{equation}\label{eq:bound}
c(g^{n})\geq 1+\frac{nk}{6(3h-1)} \ .
\end{equation}
    \end{thm}
\begin{proof}
We shall argue as in~\cite{BK}.
    As the curves $a_{i}$ are disjointly embedded in $\Sigma_{h}$, we 
    can construct a smooth Lefschetz fibration over the $2$-disk 
    $D^{2}$ with precisely one singular fiber $F_{0}$ and with vanishing 
    cycles $a_{1},\ldots,a_{k}$. The monodromy around the boundary of 
    the disk is $g$. As the Dehn twists along the $a_{i}$ are all 
    right-handed, there is a symplectic (in fact K\"ahler) structure 
    on the total space compatible with the fibration, cf.~\cite{Go}.
    Pulling back the fibration under the map $z \mapsto z^{n}$ 
    and taking the minimal resolution, we obtain a symplectic Lefschetz 
    fibration over $D^{2}$ with only one singular fiber having $k\cdot n$ 
    vanishing cycles, such that there are $n$ parallel copies of each $a_{i}$.
    The monodromy of this fibration around the boundary of the disk is 
    $g^{n}$. This can be expressed as a product of $c(g^{n})$ many commutators, 
    so that we can find a smooth surface bundle with fiber $\Sigma_{h}$ over 
    a surface of genus $c(g^{n})$ with one boundary component and the same 
    restriction to the boundary. Let $X$ be the Lefschetz fibration over 
    the closed surface $B$ of genus $c(g^{n})$ obtained by gluing together
    the two fibrations along their common boundary. 

By construction, $X$ is symplectic. As all the $a_{i}$ are 
homotopically essential, $X$ is relatively minimal, so that we 
can apply Theorem~8 of~\cite{BK}. In the notation of~\cite{BK} we 
have $D=1$ because there is only one singular fiber, and $N$, the 
total number of irreducible components of singular fibers, is $N=k\cdot 
n + c$, with $c$ a constant independent of $n$ which only depends on 
the configuration of curves $a_{1},\ldots,a_{k}\subset\Sigma_{h}$. 
Theorem~8 of~\cite{BK}
gives $k\cdot n\leq 6(3h-1)(c(g^{n})-1) +6c$. Pulling back the 
fibration to large degree coverings of the base we finally obtain 
$k\cdot n\leq 6(3h-1)(c(g^{n})-1)$ as claimed.
\end{proof}
Combining Theorem~\ref{t:main} with Proposition~\ref{p:tor}, we obtain:
\begin{cor}\label{c:Dehn}
    If $g\in\Gamma_{h}$ is the product of $k$ right-handed Dehn twists 
    along homotopically essential disjoint curves, then its stable torsion 
    length is bounded below by 
    $$
    \vert\vert g\vert\vert_{T}\geq\frac{k}{3(3h-1)} \ .
    $$
    \end{cor}
{\it A fortiori}, the torsion length is unbounded on $\Gamma_{h}$.
This answers a question raised by Brendle and Farb in their recent 
paper~\cite{BrFa}. Motivated by their paper, I had first proved the Corollary 
for powers of a single Dehn twist in a different way. Namely, it is known that 
$\Gamma_{h}$ contains only finitely many conjugacy classes of elements of 
finite 
order\footnote{This follows from a result of Harvey~\cite{H}. See~\cite{K} 
for an alternative proof.}. Therefore, the commutator length is bounded on 
torsion elements. As the commutator lengths of powers of Dehn twists grow 
linearly with the exponent by the special case of 
Theorem~\ref{t:main} proved in~\cite{EK}, the torsion lengths must also grow 
linearly\footnote{Compare also~\cite{Kork2}.}.

The existence of a non-trivial quasi-homomorphism guaranteed by the 
special case of Theorem~\ref{t:main} proved in~\cite{EK} implies that 
mapping class groups do not have the property $(TT)$ discussed by 
Monod in Section~13.2 of~\cite{Monod}, although it says nothing about 
Kazhdan's property $(T)$, of which $(TT)$ is a strengthening.

Corollary~\ref{c:Dehn} is a special case of a result which holds for 
the stable $S$-length for any $S$. Indeed, for products $g$ of commuting 
right-handed Dehn twists Theorem~\ref{t:main} implies the existence of a 
homogeneous quasi-homomorphism $\varphi$ with $\varphi(g)\neq 0$. 
This in turn implies, by the estimate~\eqref{e:qslength}, that $g$ has positive 
stable $S$-length for any $S$ on which $\varphi$ is bounded. If we 
take for $S$ a (symmetric) finite generating set of the mapping class 
group, we find that $l_{S}(g^{n})$ grows linearly with $n$. This last 
statement was proved earlier for all elements of infinite order in 
$\Gamma_{h}$ by Farb--Lubotzky--Minsky~\cite{FLM}, cf.~\cite{Kork}. In fact, 
Farb--Lubotzky--Minsky~\cite{FLM} give a direct argument estimating 
the word length only for certain reducible mapping classes, and then 
appeal to earlier work of Mosher~\cite{Mosher} for the case when $g$ 
has a power that is pseudo-Anosov, at least on a subsurface of 
$\Sigma_{h}$.

There is an important difference between the argument of~\cite{FLM} 
and ours: the former is insensitive to the chirality of Dehn twists, 
whereas the latter is not. Indeed, if we allow both right-handed Dehn twists 
and their inverses, then the proof of Theorem~\ref{t:main} breaks 
down. The corresponding fibration is no longer symplectic, and is not 
a Lefschetz fibration but an achiral Lefschetz fibration. We now show 
that this difficulty in the proof of Theorem~\ref{t:main} cannot be 
overcome by other arguments.

\begin{ex}
    Let $a\subset\Sigma_{h}$ be an essential simple closed curve, and 
    $f\colon\Sigma_{h}\rightarrow\Sigma_{h}$ a diffeomorphism\footnote{which 
    we identify with its mapping class} with $a\cap f(a)=\emptyset$. 
    Let $b=f(a)$, and $g=t_{a}t_{b}^{-1}$, where $t_{c}$ 
    denotes the right-handed Dehn twist along $c$. Then 
    $$
    g^{n}=(t_{a}t_{b}^{-1})^{n}=t_{a}^{n}t_{b}^{-n}=(f^{-1}t_{b}f)^{n}t_{b}^{-n}
    =[f^{-1},t_{b}^{n}] \ ,
    $$
    and thus the commutator length of $g^{n}$ is constant $=1$, 
    although $g$ has infinite order in $\Gamma_{h}$.
    
    If we take for $f$ an involution, then $g$ is the product of the 
    two involutions $f$ and $t_{b}ft_{b}^{-1}$. Thus, by 
    Example~\ref{ex:involutions}, the  torsion length of $g^{n}$ is 
    equal to $2$ for all $n\neq 0$. In particular, its stable torsion 
    length vanishes.
    \end{ex}
The previous example is inspired by the work of McCarthy and Papadopoulos~\cite{MP} 
on involutions in mapping class groups. They have also shown that there are 
products of involutions which are pseudo-Anosov. These are analogous to the 
Anosov diffeomorphism of the torus in Example~\ref{ex:SL}. They are elements of 
infinite order in $\Gamma_{h}$ which are conjugate to their inverses, 
and have trivial stable torsion and commutator lengths, although they 
have positive stable $S$-length for every finite $S$. This shows that 
there is no chance to find homogeneous quasi-homomorphisms which are nonzero on 
an arbitrarily fixed pseudo-Anosov mapping class. Nevertheless, in their proof 
of the following theorem, Bestvina--Fujiwara~\cite{BF} have constructed many 
quasi-homomorphisms which do not vanish on certain pseudo-Anosov elements.
\begin{thm}[\cite{BF}]\label{t:BF}
    Let $G$ be any non-virtually Abelian subgroup of a mapping class 
    group. Then the space of homogeneous quasi-homomorphisms on $G$ is 
    infinite-dimensional.
    \end{thm}
Combining this result for the case when $G$ is a mapping class group 
itself with Proposition~\ref{p:BG}, we obtain the following:
\begin{cor}\label{c:BG}
The mapping class groups $\Gamma_{h}$ are not weakly boundedly 
generated for any $h\geq 1$.
\end{cor}
Such a statement, with weak bounded generation replaced by bounded 
generation, was first proved by Farb--Lubotzky--Minsky~\cite{FLM} 
using the theory of $p$-adic analytic groups.

One can contrast the above results on the existence of homogeneous 
quasi-homomorphisms on mapping class groups with the results of 
Burger--Monod~\cite{BM}, who proved that lattices in higher rank Lie 
groups do not admit any homogeneous quasi-homomorphisms; compare 
also~\cite{Monod}. Because of their results, the special case of 
Theorem~\ref{t:main} proved in~\cite{EK} suffices to prove the 
non-arithmeticity of mapping class groups in a way which is different 
from Ivanov's proofs, except for the fact that one has to appeal to 
one of the usual arguments to exclude rank $1$ lattices, see the discussion 
in Chapter~9 of~\cite{I}. As explained by Bestvina--Fujuwara~\cite{BF},
using the existence of quasi-homomorphisms on all non-virtually Abelian 
subgroups of mapping class groups in conjunction with the results of 
Burger--Monod gives a new proof of the superrigidity result of 
Farb--Masur~\cite{FM}, which is independent of~\cite{FM} and~\cite{KM}. 
This argument does not require the full strength of Theorem~\ref{t:BF}, 
as one only needs the non-triviality (rather than the infinite-dimensionality) 
of the space of homogeneous quasi-homomorphisms on every non-virtually 
Abelian subgroup of a mapping class group. In fact, the above proof of 
Corollary~\ref{c:BG} is the only application I know of the 
infinite-dimensionality of the space of homogeneous quasi-homomorphisms 
on $\Gamma_{h}$.

\begin{rem}
    The discussion in this section applies to mapping class groups of 
    surfaces with marked points and/or boundaries, but it does not 
    apply to extended mapping class groups consisting of isotopy 
    classes of not necessarily orientation-preserving diffeomorphisms. 
    In an extended mapping class group a Dehn twist is conjugate to 
    its inverse, because it is a product of two involutions, see 
    A'Campo~\cite{NAC}. Thus, every quasi-homomorphism defined on an 
    extended mapping class group must vanish on a Dehn twist, 
    cf.~Example~\ref{ex:involutions}.
    \end{rem}

\medskip\noindent
{\bf Acknowledgments:} I am grateful to B.~Farb for pointing out the 
question in~\cite{BrFa} about the torsion lengths of elements of mapping 
class groups. Having answered this question in an ad hoc way, I was 
spurred to think systematically about such questions by the 
appearance of M.~Korkmaz's preprint~\cite{Kork2}. Thanks also to 
M.~Burger for useful email communications.

\bibliographystyle{amsplain}

\bigskip

\end{document}